\documentclass[11pt]{article}

\usepackage{fullwidth}
\usepackage{color}
\usepackage[color,matrix,arrow]{xy}
\usepackage[top=1.5in, bottom=1.5in, left=1in, right=1in]{geometry}
\usepackage{amsgen}
\usepackage{amsmath}
\usepackage{amstext}
\usepackage{amsbsy}
\usepackage{amsopn}
\usepackage{amsfonts}
\usepackage{amssymb}
\usepackage{eepic}
\usepackage{graphicx}
\usepackage{epsf}
\usepackage{pstricks}
\xyoption{all}

\def\Box{\square}

\def\tra#1{\smash{\mathop{\mid\kern
-1pt\joinrel\relbar\joinrel\relbar}\limits^{*}_{#1}}}
\def\longtra#1{\smash{\mathop{\mid\kern
-1pt\joinrel\relbar\joinrel\relbar\joinrel\relbar}\limits^{*}_{#1}}}
\def\vlongtra#1{\smash{\mathop{\mid\kern
-1pt\joinrel\relbar\joinrel\relbar\joinrel\relbar\joinrel\relbar}\limits^{*}_{#1}}}
\def\vvlongtra#1{\smash{\mathop{\mid\kern
-1pt\joinrel\relbar\joinrel\relbar\joinrel\relbar\joinrel\relbar\joinrel\relbar}\limits^{*}_{#1}}}
\def\vvvlongtra#1{\smash{\mathop{\mid\kern
-1pt\joinrel\relbar\joinrel\relbar\joinrel\relbar\joinrel\relbar\joinrel\relbar\joinrel\relbar}\limits^{*}_{#1}}}
\def\etra#1{\smash{\mathop{\mid\kern
-1pt\joinrel\relbar\joinrel\relbar}\limits_{#1}}}

\def\A{{\cal{A}}}

\def\iff{\Leftrightarrow}
\def\Rw{\Rightarrow}
\def\oo{\overline}
\def\wt{\widetilde}

\def\N{\mathbb{N}}

\def\rat{\mbox{Rat}}

\def\rec{\mbox{Rec}}

\def\mt{\mbox{MT}}
\def\max{\mbox{max}}
\def\min{\mbox{min}}

\def\Z{\mathbb{Z}}

\def\p{\varphi}

\def\inv{^{-1}}
\def\la{\langle}
\def\ra{\rangle}

\def\bi{\begin{itemize}}
\def\ei{\end{itemize}}
\def\beq{\begin{equation}}
\def\eeq{\end{equation}}

\def\J{{\cal{J}}}

\def\xr{\xrightarrow}
\def\ds{\displaystyle}

\newtheorem{T}{Theorem}[section]
\newcommand{\bt}{\begin{T}}
\newcommand{\et}{\end{T}}
\newcommand{\ftd}{$\square$\end{T}}

\newtheorem{Proposition}[T]{Proposition}
\newcommand{\bp}{\begin{Proposition}}
\newcommand{\ep}{\end{Proposition}}
\newcommand{\fpd}{$\square$\end{Proposition}}

\newtheorem{Lemma}[T]{Lemma}
\newcommand{\bl}{\begin{Lemma}}
\newcommand{\el}{\end{Lemma}}
\newcommand{\fld}{$\square$\end{Lemma}}

\newtheorem{Corol}[T]{Corollary}
\newcommand{\bc}{\begin{Corol}}
\newcommand{\ec}{\end{Corol}}
\newcommand{\fcd}{$\square$\end{Corol}}

\newtheorem{Result}[T]{Result}
\newcommand{\br}{\begin{Result}}
\newcommand{\er}{\end{Result}}
\newcommand{\frd}{$\square$\end{Result}}

\newtheorem{Example}[T]{Example}
\newcommand{\be}{\begin{Example}}
\newcommand{\ee}{\end{Example}}

\newtheorem{Problem}[T]{Problem}
\newcommand{\bq}{\begin{Problem}}
\newcommand{\eq}{\end{Problem}}

\newtheorem{Remark}[T]{Remark}
\newcommand{\brm}{\begin{Remark}}
\newcommand{\erm}{\end{Remark}}

\newcommand{\proof}
   {\par\medbreak\noindent{\bf Proof}.\enspace}

\newcommand{\qed}{\hfill
$\Box$
\par\bigbreak}

\normalbaselines
\def\abstract#1{\par\bigskip
\begingroup\small
\baselineskip=12truept
\begin{center}ABSTRACT\end{center}
\par\medskip\par\noindent
\null\hfill\hbox{\vbox{\hsize=5truein\noindent#1}}
\hfill\null\par\endgroup\par}

\title{On the rational subsets of the monogenic free inverse monoid}
\author{{\bf Pedro V. Silva}\\ $ $\\
{\em Centro de
Matem\'{a}tica, Faculdade de Ci\^{e}ncias, Universidade do
Porto,}\\ {\em R. Campo Alegre 687, 4169-007 Porto, Portugal}\\
{\em email:} pvsilva@fc.up.pt}

\date{\today}

\begin{document}
\maketitle

\begin{center}\small

\bigskip

\end{center}

\abstract{We prove that the equality problem is decidable for rational subsets of the monogenic free inverse monoid $F$. It is also decidable whether or not a rational subset of $F$ is recognizable. We prove that a submonoid of $F$ is rational if and only if it is finitely generated. We also prove that the membership problem for rational subsets of a finite $\J$-above monoid is decidable, covering the case of free inverse monoids.
}

\section{Introduction}

Inverse semigroups have a special place in semigroup theory: if semigroups are represented through transformations of a set and groups through permutations, inverse semigroups are represented through partial permutations of a set (the Vagner-Preston representation). This simple fact makes them almost ubiquituous in Mathematics. In this paper, we work with inverse monoids for commodity, but all the results hold for inverse semigroups as well (we only need to remove the identity element).

Inverse monoids constitute a variety of algebras of signature $(2,1,0)$, so free inverse monoids are bound to exist. We denote by $FIM_A$ the free inverse monoid on the alphabet $A$.
However, normal forms (and consequently a solution of the word problem) were only obtained in the early seventies, through the independent work of Munn and Scheiblich \cite{Mun,Sch}. We will favour in this text the elegant solution provided by Munn, using a special type of finite automata known as Munn trees.

The monogenic case was settled earlier by Gluskin in 1957 \cite{Glu}. The monogenic free inverse monoid, which will be denoted by $F$, is a noncommutative monoid which can be quite challenging in comparison with other monogenic free objects such as $\N$ or $\Z$. It is the most studied monogenic inverse monoid alongside with the famous {\em bicyclic monoid}, introduced by Lyapin in 1953 \cite{Lya}. The latter can be described through the inverse monoid presentation $\la a \mid aa\inv = 1\ra$, or as the transformation monoid generated by the shift map $$\begin{array}{rcl}
\sigma:\N&\to&\N\setminus \{ 0\}\\
n&\mapsto&n+1
\end{array}$$
and its inverse. It also plays an important role in ring theory \cite{PS}.

Which other decidability results can we consider for free inverse monoids? In \cite{OS}, Oliveira and the author proved that every finitely generated inverse submonoid of $F$ is finitely presented. As a consequence, the isomorphism problem for $F$ is decidable.

But the word problem can be generalized in a different direction using the concept of rational subset. Rational subsets of a monoid can be viewed as finitely generated subsets in some precise sense, admitting also a description in terms of finite automata. With respect to connections between free inverse monoids and automata theory, see \cite{Sil}.

With increasing level of difficulty, classical problems are the membership problem and the equality problem for rational subsets. The equality problem is in general a tough problem, often undecidable. As far as free inverse monoids are concerned, the membership problem is relatively easy to settle for arbitrary alphabets, the equality problem is harder. But we can solve it in the monogenic case.

The main obstacle here is that there is no obvious normal form for describing rational subsets of $F$. In the case of rational languages, we have minimal automata, and the bicyclic case follows from a generalization of Benois's Theorem due to S\'enizergues \cite{Sen}, which outputs also a minimal finite automaton. But no such normal form is known for rational subsets of $F$.

Still, we succeed on proving decidability theorems for $F$. The membership problem for rational subsets follows from a general result on finite $\J$-above monoids, proved in Section 3. Then in Section 4 we prove the technical {\em cut and paste lemma} as preparation for the decision algorithm for inclusion in Section 5, which solves also the equality problem. 

Given a finite generated monoid $M$, its recognizable subsets are rational, but the converse does not hold in general. This is an important concept in language theory because it collects those subsets which can be recognized through a finite monoid. In Section 6 we show that it is decidable whether or not a rational subset of $F$ is recognizable. 

Anisimov and Seifert proved in 1975 that a subgroup of a free group is rational if and only if it is finitely generated \cite{AS}. We generalize their theorem in Section 7 proving that a submonoid of $F$ is rational if and only if it is finitely generated. Finally, we propose in Section 8 a couple of open problems.

\section{Preliminaries}

We assume that $0 \in \N$. For every $n \geq 1$, we write $[n] = \{ 1,\ldots,n\}$.

The reader is assumed to have basic knowledge of inverse semigroup theory and automata theory, being respectively referred to \cite{Law} and \cite{Ber,Sak} for that purpose.

\subsection{Rational subsets}

Given a finite alphabet $A$, an $\A$-{\em automaton} is a
structure of the form ${\cal{A}} = (Q,I,T,E)$, where:
\bi
\item
$Q$ is the set of vertices,
\item
$I,T \subseteq Q$ are the subsets of initial and terminal vertices, respectively,
\item
$E \subseteq Q \times A \times Q$ is the set of edges.
\ei
This automaton is finite if $Q$ is finite.
Its language is defined as
$$L(\A) = \{ u \in A^* \mid \mbox{ there exists some path $I \ni i \xr{u} t \in T$ in }\A\}.$$
In pictures, we will identify initial vertices by an incoming arrow $\rightarrow \cdot$ or $\circ$, and terminal vertices by
an outgoing arrow $\cdot \rightarrow$ or $\bullet$.

We say that $L \subseteq A^*$ is a {\em rational language} if $L$ can be obtained from finite subsets of $A^*$ using finitely many times the operators union, product and star:
$$(X,Y) \mapsto X \cup Y,\quad (X,Y) \mapsto XY,\quad X \mapsto X^*,$$
where $X^* = \ds\bigcup_{n \geq 0} X^n$. 

The classical Kleene's Theorem states that an $A$-language $L$ is rational if and only if $L = L(\A)$ for some finite $A$-automaton $\A$. 

Given a monoid $M$, we say that $K \subseteq M$ is a {\em rational subset of} $M$ if $K$ can be obtained from finite subsets of $M$ using finitely many times the operators union, product and star. We denote by $\rat(M)$ the set of all rational subsets of $M$. A {\em rational submonoid} of $M$ is a rational subset which is also a submonoid of $M$. 

Fix an alphabet $A$ and a surjective homomorphism $\theta:A^* \to M$. Then $K \subseteq M$ is a rational subset of $M$ if and only if $K = L\theta$ for some rational language $L \subseteq A^*$. We consider the following classical decidability problems for $M$, of increasing difficulty:
\bi
\item
{\bf Word problem:} Is there an algorithm which decides, on input $u,v \in A^*$, whether or not $u\theta = v\theta$?
\item
{\bf Membership problem for rational subsets:} Is there an algorithm which decides, on input $u\in A^*$ and $L \in \rat(A^*)$, whether or not $u\theta \in L\theta$?
\item
{\bf Equality problem for rational subsets:} Is there an algorithm which decides, on input $K,L \in \rat(A^*)$, whether or not $K\theta = L\theta$?
\ei
If $M$ is finitely generated, the alphabet $A$ can be chosen to be finite. Then decidability of any of these problems does not depend on the finite alphabet $A$ or the surjective homomorphism $\theta$.

\subsection{Free inverse monoids}

Let $A$ be a nonempty alphabet. We extend $^{-1}:A \to A^{-1}: a \mapsto a^{-1}$ to an involution on the free monoid 
$\widetilde{A}^*$ through 
$$(a^{-1})^{-1} = a,\quad 1\inv = 1,\quad (ua)^{-1} = a^{-1}u^{-1}\hspace{1cm} (u \in \wt{A}^*;\;
a\in \widetilde{A})\, .$$
The {\em free inverse monoid on} $A$ is the quotient $FIM_A = \wt{A}^*/\nu$,
where $\nu$ is the congruence on $\wt{A}^*$ generated by the relation
$$\{ (ww\inv w,w) \mid w \in \wt{A}^*\} \cup
\{ (uu\inv vv\inv,vv\inv uu\inv) \mid u,v \in \wt{A}^*\}.$$
known as the {\em Wagner congruence} on $\wt{A}^*$. 

W. D. Munn provided in \cite{Mun} an elegant solution to the word problem for $FIM_A$ using automata (see also \cite{Sch} by Scheiblich), which we now describe.

Given $w = a_1\ldots a_n \in \wt{A}^*$ $(a_i \in \wt{A})$, let Lin$(w)$ denote the {\em linear automaton} of $w$:
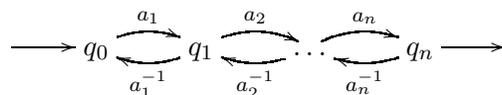
$$\xymatrix{
\ar[r] & q_0 \ar@/^/[r]^{a_1} & q_1 \ar@/^/[r]^{a_2} \ar@/^/[l]^{a_1\inv}&\ldots\ar@/^/[l]^{a_2\inv} \ar@/^/[r]^{a_n} & q_n \ar[r] \ar@/^/[l]^{a_n\inv}&
}$$

Departing from Lin$(w)$, Munn proceeds to successively identifying distinct edges of the form 
$$p \xleftarrow{a} q \xr{a} r$$
for all $a \in \wt{A}$. Since the number of edges decreases at each identification, this procedure will eventually terminate. Munn proved that the final automaton is independent of the sequence of identifications. 

Ten years later, after Stallings seminal paper \cite{Sta}, this folding operation would become known in combinatorial group theory as {\em Stallings foldings}. 

The automaton obtained by folding Lin$(w)$ became known as the {\em Munn tree} of $w$, and it will be denoted by $\mt(w)$. It can be viewed as a legitimate tree if we view each pair of edges
$$\xymatrix{
p \ar@/^/[r]^{a} & q \ar@/^/[l]^{a\inv}
}$$
as a single undirected edge. We denote by ${\rm MT}^{{\rm o}}(w)$ the underlying graph of ${\rm MT}(w)$.

In the graphical representation of both Lin$(w)$ and $\mt(w)$, it suffices to depict those edges with label in $A$:

\be
\label{mtree}
Let $w = ab^2b\inv a\inv a^2a\inv$.
Then {\rm Lin}$(w)$ is the automaton
$$\xymatrix{
\circ \ar[r]^a & \cdot  \ar[r]^b & \cdot  \ar[r]^b & \cdot &\cdot \ar[l]_{b} &\cdot \ar[l]_{a} \ar[r]^a & \cdot \ar[r]^a &\cdot & \bullet \ar[l]_{a} 
}$$
and ${\rm MT}(w)$ is the automaton
$$\xymatrix{
&&\cdot\ar[d]^a&\\
\circ\ar[r]^a&\cdot\ar[r]^b&\bullet\ar[d]^a \ar[r]^b&\cdot\\
&&\cdot&
}$$
\ee

Munn trees provide the following solution for the word problem of $FIM_A$:
 
\bt
\label{wpfree}
{\rm \cite{Mun}}
For all $u,v \in \wt{A}^*$, the following conditions are equivalent:
\bi
\item[(i)] $u\nu = v\nu$;
\item[(ii)] ${\rm MT}(u) \cong {\rm MT}(v)$.
\ei
\et

Given a monoid $M$, the $\J$-{\em order} on $M$ is defined by
$$u \leq_{\J} v \mbox{ if } u \in MvM$$
and we write $u\,\J\,v$ if $u\leq_{\J} v$ and $v \leq_{\J} u$. 

\bp
\label{jfim}
{\rm \cite{Mun}}
Let $u,v \in \wt{A}^*$. Then:
\bi
\item[(i)]
$u\nu \leq_{\J} v\nu$ if and only if ${\rm MT}^{{\rm o}}(v)$ embeds into ${\rm MT}^{{\rm o}}(u)$;
\item[(ii)]
$u\nu\,\J\, v\nu$ if and only if ${\rm MT}^{{\rm o}}(u)\cong {\rm MT}^{{\rm o}}(v)$.
\ei
\ep

We consider now the monogenic case. We denote by $F$ the free inverse monoid on the set $\{ a\}$ and let $\pi:F \to \Z$ be the homomorphism defined by $a\pi = 1$. 
The underlying graph of a Munn tree in $F$ is necessarily of  the form
$$\xymatrix{
\cdot \ar[r]^a & \cdot \ar[r]^a & \cdot \ar[r]^a & \ldots \ar[r]^a & \cdot \ar[r]^a &\cdot
}$$
and then we must specify the initial and the terminal vertex (if they coincide, we have an idempotent). It follows easily that
every $u \in F$ admits a unique factorization of the form $u = a^{-m}a^{m+k}a^{-k}a^{u\pi}$ with $m,k\in \N$ and $-m \leq u\pi \leq k$. We shall write $u\lambda = -m$ and $u\rho = k$. It is routine to check that
\beq
\label{prod}
(uv)\lambda = \min\{ u\lambda, u\pi+v\lambda\}, \quad (uv)\rho = \max\{ u\rho, u\pi + v\rho\},\quad (uv)\pi = u\pi + v\pi
\eeq
holds for all $u,v \in F$. We shall abuse notation by writing $w\xi = w\theta\xi$ for all $w \in A^*$ and $\xi \in \{ \lambda,\rho,\pi\}$.

The {\em norm} of $w \in F$ is defined as $||w|| = w\rho-w\lambda$. This is the number of positive edges in $\mt(w)$. It is clear that $||uv|| \geq ||u||,||v||$ holds for all $u,v \in F$.

\section{Finite $\J$-above monoids}

For every $u \in M$, let $u\gamma = \{ v \in M \mid u \leq_{\J} v\}$. We say that $M$ is {\em finite $\J$-above} if $u\gamma$ is finite for every $u \in M$.

\bt
\label{fja}
Let $M$ be a finitely generated, finite $\J$-above monoid such that:
\bi
\item[(i)] $M$ has decidable word problem;
\item[(ii)] the $\J$-order on $M$ is decidable.
\ei
Then the membership problem for rational subsets of $M$ is decidable.
\et

\proof
Since $M$ is finitely generated, there exists some finite alphabet $A$ and some surjective homomorphism $\theta:A^* \to M$. Let $\A = (Q,I,T,E)$ be a finite $A$-automaton and let $u \in M$. We present an algorithm to decide whether or not $u \in (L(\A))\theta$. Since 
$$L(\A) = \bigcup_{q_0\in I}\bigcup_{t \in T} L(Q,q_0,t,E),$$
we may assume that $I = \{ q_0\}$ and $T = \{ t\}$.

We build a finite sequence $(r_0,u_0),\ldots,(r_m,u_m)$ on $Q \times u\gamma$ under the following rules:
\bi
\item
$(r_0,u_0) = (q_0,1)$;
\item
if $(r_0,u_0),\ldots,(r_i,u_i)$ are defined, and there exist $j \in \{ 0,\ldots,i\}$, $a \in A$ and $q \in Q$ such that
\bi
\item
$(r_j,a,q) \in E$, 
\item
$u_j(a\theta) \in u\gamma$,
\item
$(q,u_j(a\theta)) \notin \{ (r_0,u_0),\ldots,(r_i,u_i) \}$, 
\ei
we choose $(r_{j+1},u_{j+1}) = (q,u_j(a\theta))$.
\ei
Since $Q$ and $u\gamma$ are both finite and all the terms of the sequence are distinct, every such sequence (which we do not claim to be unique) is necessarily finite (so the procedure described above is bound to terminate). In view of our decidability conditions, we can compute such a sequence. Moreover, we prove by induction that, for $i = 0,\ldots,m$, there exist some $s \geq 0$, $0 = k_0 < k_1 < \ldots < k_s = i$ and a path
$$r_{k_0} \xr{a_1} r_{k_1} \xr{a_2} \ldots \xr{a_s} r_{k_s}$$
in $\A$ such that $(a_1\ldots a_j)\theta = u_{k_j}$ for $j = 0,\ldots,s$.

Indeed, the claim holds trivially for $i = 0$. Let $i \in [m]$ and assume it holds for $0,\ldots,i-1$.  It follows from the construction of the sequence that there exist some $\ell < i$ and some $a \in A$ such that $(r_{\ell},a,r_i) \in E$ and $u_{\ell}(a\theta) = u_i$. By the induction hypothesis, there exist some $s \geq 0$, $0 = k_0 < k_1 < \ldots < k_s = \ell$ and a path
$$r_{k_0} \xr{a_1} r_{k_1} \xr{a_2} \ldots \xr{a_s} r_{k_s}$$
in $\A$ such that $(a_1\ldots a_j)\theta = u_{k_j}$ for $j = 0,\ldots,s$. Now we consider $0 = k_0 < k_1 < \ldots < k_s = \ell < i$
and the path
$$r_{k_0} \xr{a_1} r_{k_1} \xr{a_2} \ldots \xr{a_s} r_{k_s} = r_{\ell} \xr{a} r_i$$
in $\A$. We have $(a_1\ldots a_j)\theta = u_{k_j}$ for $j = 0,\ldots,s$. On the other hand, 
$$(a_1\ldots a_sa)\theta = u_{k_s}(a\theta) = u_{\ell}(a\theta) = u_i,$$
hence our induction is complete. 

Now we show that
\beq
\label{fja1}
u \in (L(\A))\theta\quad\mbox{if and only if}\quad
(t,u) = (r_i,u_i)\mbox{ for some }i \in \{ 0,\ldots,m\}.
\eeq

Suppose first that $u \in (L(\A))\theta$. Then $u = v\theta$ for some $v\in L(\A)$. 
Write $v = a_1\ldots a_k$ with $a_1,\ldots,a_k \in A$. Then there exists some path
$$q_0 \xr{a_1} q_1 \xr{a_2} \ldots \xr{a_k} q_k = t$$
in $\A$. We prove by induction that
\beq
\label{fja2}
\forall i \in \{ 0,\ldots,k\}\, \exists j \in \{0,\ldots,m\}:
(q_i,(a_1\ldots a_i)\theta) = (r_j,u_j).
\eeq

This holds trivially for $i = 0$. Assume now that $i \in [k]$ and the claim holds for $i-1$. Then  
$(q_{i-1},(a_1\ldots a_{i-1})\theta) = (r_j,u_j)$ for some $j \in \{0,\ldots,m\}$. Since 
$(a_1\ldots a_i)\theta \geq_{\J} v\theta = u$, we have $(a_1\ldots a_i)\theta \in u\gamma$. Now it follows from the definition of the sequence $(r_0,u_0),\ldots,(r_m,u_m)$ that $(q_i,(a_1\ldots a_i)\theta) = (r_k,u_k)$ for some $k \in \{0,\ldots,m\}$ and so (\ref{fja2}) holds.

In particular, 
$(t,u) = (q_k,(a_1\ldots a_k)\theta) = (r_s,u_s)$ for some $s \in \{0,\ldots,m\}$. Therefore the direct implication of (\ref{fja1}) holds.

Conversely, suppose that $(t,u) = (r_i,u_i)$ for some $i \in \{ 0,\ldots,m\}$. By our previous claim, there exist some $s \geq 0$, $0 = k_0 < k_1 < \ldots < k_s = i$ and a path
$$q_0 = r_0 = r_{k_0} \xr{a_1} r_{k_1} \xr{a_2} \ldots \xr{a_s} r_{k_s} = r_i = t$$
in $\A$ such that $(a_1\ldots a_j)\theta = u_{k_j}$ for $j = 0,\ldots,s$. In particular, 
$$u = u_i = u_{k_s} = (a_1\ldots a_s)\theta \in (L(\A))\theta$$
and so (\ref{fja1}) holds.

Since the sequence $(r_0,u_0),\ldots,(r_m,u_m)$ is computable and $M$ has decidable word problem, we can decide whether or not $(t,u) = (r_i,u_i)$ for some $i \in \{ 0,\ldots,m\}$. By (\ref{fja1}), we can decide whether or not $u \in (L(\A))\theta$. Therefore the membership problem for rational subsets of $M$ is decidable.
\qed

\bc
\label{mpr}
The membership problem for rational subsets of a free inverse monoid is decidable.
\ec

\proof
Let $L$ be a rational subset of $FIM_A$. Only finitely many letters of $A$ occur in $L$, hence we may assume that $A$ is a finite alphabet. By Theorem \ref{wpfree}, $FIM_A$ has decidable word problem. And it follows easily from Proposition \ref{jfim} that
$FIM_A$ is finite $\J$ above and the $\J$-order on $FIM_A$ is decidable.
Therefore the claim follows from Theorem \ref{fja}.
\qed

\section{Cut and paste}

The main goal of this section is to prove a helpful {\em cut and paste lemma}.

From now on, we write $A = \{ a,a\inv\}$ and let $\theta:A^* \to F$ be the canonical homomorphism.
The next result is a basic tool which will be used implicitly throughout the paper without further reference.

\bl
\label{prefix}
Let $v \in A^*$ and let $p$ be an integer such that $v\theta\lambda \leq p \leq v\theta\rho$. Then there exists some prefix $w$ of $v$ such that $w\theta\pi = p$.
\el

\proof
Write $v = v_1\ldots v_k$ with $v_1,\ldots,v_k \in A$. Let $w_i = v_1\ldots v_i$ for $i = 0,\ldots,k$. If $p = 0$, then 
$w_0\theta\pi = p$. Suppose now that $p > 0$. 
Since $|w_i\theta\rho - w_{i-1}\theta\rho| \leq 1$ for every $i$, and $p \leq w_k\theta\rho$, there exists some $j$ such that $w_j\theta\rho = p$. Assume that $j$ is minimum. Then $w_j\theta\rho = w_{j-1}\theta\rho +1$ and this can only happen if $w_j\theta\pi = p$. The case $p < 0$ is treated similarly, replacing $\rho$ by $\lambda$.
\qed

We introduce the automorphism $\alpha:F \to F$ which permutes $a$ and $a\inv$, and $\beta:F \to F$ defined by $w\beta = (w\alpha)\inv$. This is an anti-automorphism of order 2, since $\beta^2 = 1$ and $(uv)\beta = (v\beta)(u\beta)$ for all $u,v \in F$. We can describe the action of $\beta$ in terms of Munn trees, distinguishing the cases $u\pi \geq 0$ and $u\pi \leq 0$:
$$\xymatrix{
\mt(u):& \ar[rr]^{a^i} && \circ \ar[r]^{a^j} & \bullet \ar[rrr]^{a^k} &&&\\
\mt(u\beta):&&& \bullet \ar[ll]_{a^i} & \circ \ar[l]_{a^j} &&& \ar[lll]_{a^k}\\
\mt(v):& \ar[rr]^{a^i} && \bullet \ar[r]^{a^j} & \circ \ar[rrr]^{a^k} &&&\\
\mt(v\beta):&&& \circ \ar[ll]_{a^i} & \bullet \ar[l]_{a^j} &&& \ar[lll]_{a^k}
}$$
 Now it follows easily that, for every $u \in F$ such that $u\pi \geq 0$, we have
 \beq
 \label{beta1}
 u\beta\lambda = -k = -(u\rho-u\pi),\quad u\beta\rho-u\beta\pi = i = -u\lambda,\quad u\beta\pi = j = u\pi.
 \eeq
 Similarly, if $v\pi \leq 0$, we have
 \beq
 \label{beta2}
v\beta\lambda -v\beta\pi = -k = -v\rho, \quad v\beta\rho = i = -(v\lambda-v\pi),\quad v\beta\pi = -j = v\pi.
\eeq
 Note that we can easily produce a finite $A$-automaton recognizing $L\beta$ by reversing the arrows in $\A$ and exchanging the initial with the terminal vertices.
 
For every $n \geq 1$, write
\bi
\item
$G_{n,1} = \{ u \in F \mid u\pi \geq 0,\; u\lambda \leq -n\}$,
\item
$G_{n,2} = \{ u \in F \mid u\pi \geq n\}$,
\item
$G_{n,3} = \{ u \in F \mid u\pi \geq 0,\; u\rho-u\pi \geq n\}$.
\ei
Let $i \in [3]$. We define a function $\xi_{n,i}:G_{n,i} \to F$ as follows. Let $u \in G_{n,i}$ and write $u' = u\xi_{n,i}$.
\bi
\item
If $i = 1$, then 
$u'\lambda = u\lambda+n$, $u'\rho = u\rho$ and $u'\pi = u\pi$:
$$\xymatrix{
\mt(u): & \ \ar[rrrr]^{a^i} &&&& \circ \ar[r]^{a^j} & \bullet \ar[r]^{a^k} &\\
\mt(u'): &&& \ \ar[rr]^{a^{i-n}} && \circ \ar[r]^{a^j} & \bullet \ar[r]^{a^k} &
}$$
\item
If $i = 2$, then 
$u'\lambda = u\lambda$, $u'\rho = u\rho-n$ and $u'\pi = u\pi-n$:
$$\xymatrix{
\mt(u): & \ \ar[r]^{a^i} & \circ \ar[rrrr]^{a^j} &&&& \bullet \ar[r]^{a^k} &\\
\mt(u'): & \ \ar[r]^{a^i} & \circ \ar[rr]^{a^{j-n}} && \bullet \ar[r]^{a^k} &&&
}$$
\item
If $i = 3$, then $u'\lambda = u\lambda$, $u'\rho = u\rho-n$ and $u'\pi = u\pi$:
$$\xymatrix{
\mt(u): & \ \ar[rr]^{a^i} && \circ \ar[r]^{a^j} & \bullet \ar[rrrr]^{a^k} &&&&\\
\mt(u'): & \ \ar[rr]^{a^i} && \circ \ar[r]^{a^j} & \bullet \ar[rr]^{a^{k-n}} &&&&
}$$
\ei

Note that
\beq
\label{beta4}
\xi_{n,1}\beta = \beta\xi_{n,3}\quad\mbox{and}\quad \xi_{n,3}\beta = \beta\xi_{n,1}.
\eeq

We prove next the cut and paste lemma:

\bl
\label{xix}
Let $L \in {\rm Rat}(F)$. Then there exist some computable $n' \geq n \geq 1$ such that 
\beq
\label{xix1}
u \in L \iff u\xi_{n,i} \in L
\eeq
holds for all $i \in [3]$ and $u \in G_{n',i}$.
\el

\proof
Assume that $L = (L(\A))\theta$ for the finite $A$-automaton $\A = (Q,I,T,E)$. Write $m = |Q|$. We consider three cases.
We can always identify the case we are in because $L\pi \in \rat(\Z)$ and by Benois' Theorem we can always build a finite automaton recognizing the reduced forms of $L\pi$ \cite[Theorem 4.3]{BS1}.

\medskip

\noindent
\underline{Case 1}: $L\pi \cap \N$ is finite.

\medskip

Let $K = \{ u \in L \mid u\pi \geq 0\}$. Let $n = 2m$ and $n' = 4m$.

Let $u \in K$ and $x \in u\theta\inv \cap L(\A)$.
Suppose that $||u|| \geq n = 2m$. Then there exists a path in $\A$ of the form
\beq
\label{afo1}
I \ni q_0 \xr{x'} p_0 \xr{x_1} p_1 \xr{x_2} \ldots \xr{x_m} p_{m} \xr{x''} t \in T
\eeq
with $x = x'x_1\ldots x_{m}x''$ and $x_j\pi = 1$ for $j \in [m]$. By the Pigeonhole Principle, there exist $0 \leq j_1 < j_2 \leq m$ such that $p_{j_1} = p_{j_2}$. Now 
$$y_k = x'x_1\ldots x_{j_1}(x_{j_1+1}\ldots x_{j_2})^k x_{j_2+1}\ldots x_{m}x'' \in L(\A)$$
for every $k \geq 1$. Moreover, $y_k\pi = u\pi + (k-1)(j_2-j_1)$, a contradiction in Case 1. Therefore $||u|| < n$.

Let $u \in G_{2n,i}$. Then $u \in L$ implies $u \in K$ and so $||u|| < n$, a contradiction. Similarly, 
$u\xi_{n,i} \in L$ implies $u\xi_{n,i} \in K$, hence $||u\xi_{n,i}|| < n$  and so $||u|| < 2n$, also a contradiction. Therefore (\ref{xix1}) holds in this case.

\medskip

\noindent
\underline{Case 2}: $L\pi \cap (-\N)$ is finite.

\medskip

Let $n = m!$ and $n' = (m+3)m!$. 
Adapting the argument used in Case 1, we see that every $u \in L$ must satisfy $u\lambda > -m$ and $u\rho-u\pi < m$. Hence 
$L \cap (G_{n',i} \cup G_{n',i}\xi_{n,i}) = \emptyset$ for $i = 1,3$ and so (\ref{xix1}) holds for $i = 1,3$.
Thus we only have to deal with the case $i = 2$.

We start showing that
\beq
\label{afo2}
\mbox{if $u\xi_{n,2} \in G_{3n,2} \cap L$, then $u \in L$.}
\eeq

Indeed, take $x \in u\xi_{n,2}\theta\inv \cap L(\A)$.
Then there exists a path in $\A$ of the form (\ref{afo1})
with $x = x'x_1\ldots x_{m}x''$, $x'\pi = m$ and $x_j\pi = 1$ for $j \in [m]$. Hence there exist $0 \leq j_1 < j_2 \leq m$ such that $p_{j_1} = p_{j_2}$. Let $k = \frac{n}{j_2-j_1}$. Then 
$$y = x'x_1\ldots x_{j_1}(x_{j_1+1}\ldots x_{j_2})^{k+1} x_{j_2+1}\ldots x_{m}x'' \in L(\A).$$

Clearly, $y\pi = x\pi + k(j_2-j_1) = x\pi + n$. We remark that no factor $w$ of a word in $L(\A)$ can satisfy $w\pi \leq -m$, otherwise by inserting loops as we did before we would end up contradicting Case 2. Hence
$x'\pi = m$ implies that $y\lambda = x'\lambda = u\lambda$. Since 
$$x''\pi = x\pi - (x'x_1\ldots x_{m})\pi  = u\xi_{n,2}\pi -2m \geq 3n -2m \geq m,$$
a similar argument yields
$y\rho -y\pi = x''\rho-x''\pi = u\rho-u\pi$. Hence $y\theta = u$ and so $u \in L$. Therefore (\ref{afo2}) holds.

Now let $u \in G_{n',2} \subseteq G_{3n,2}$. If $u\xi_{n,2} \in L$, then $u \in L$ by (\ref{afo2}). 

Conversely, assume that $u \in L$. Take $x \in u\xi_{n,2}\theta\inv \cap L(\A)$.
Then there exists a path in $\A$ of the form (\ref{afo1})
with $x = x'x_1\ldots x_{m}x''$, $x'\pi = m$ and $x_j\pi = m!$ for $j \in [m]$. Hence there exist $0 \leq j_1 < j_2 \leq m$ such that $p_{j_1} = p_{j_2}$. Write $k = j_2-j_1$. Then 
$y = x'x_1\ldots x_{j_1}x_{j_2+1}\ldots x_{m}x'' \in L(\A)$ and as in the direct implication we get $y\lambda = x'\lambda = u\lambda$ and 
$y\rho -y\pi = x''\rho-x''\pi = u\rho-u\pi$. On the other hand, 
$$y\pi = u\pi -km! \geq (m+3)n - nm = 3n.$$
Hence
$u\xi_{n,2}^k = y\theta \in L$. By applying (\ref{afo2}) $k-1$ times, we get $u\xi_{n,2} \in L$. Therefore (\ref{xix1}) holds in this case.

\medskip

\noindent
\underline{Case 3}: $L\pi \cap \N$ and $L\pi \cap (-\N)$ are both infinite.

\medskip

Let 
$$P^+ = \{ p \in Q \mid \mbox{ there exists some loop 
$\xymatrix{
p\ar@(ur,r)^{w_p}
}$ in $\A$ with $w_p\pi > 0\}$},$$
$$P^- = \{ p \in Q \mid \mbox{ there exists some loop 
$\xymatrix{
p\ar@(ur,r)^{w_p}
}$ in $\A$ with $w_p\pi < 0\}.$}$$

Let $p \in Q$. Then 
$$p \in P^+\quad\mbox{if and only if}\quad
(L(Q,p,p,E))\pi \cap (\N \setminus \{ 0\}) \neq \emptyset.$$
Since $\rat(\Z)$ is closed under intersection \cite[Corollary 4.4]{BS1} and emptiness is decidable, we can decide whether or not this intersection is nonempty, and in that case we can compute one of its elements. Therefore we can compute $P^+$ and the corresponding $w_p$. The same holds for $P^-$.  

Let $n \geq 1$ be the least common multiple of all the $w_p\pi$ $(p \in P^+ \cup P^-)$.
Then:
\bi
\item
$\mbox{for every $p \in P^+$ there exists some loop
$\xymatrix{
p\ar@(ur,r)^{z_p}
}$ in $\A$ with $z_p\pi = n$,}$
\item
$\mbox{for every $p \in P^-$ there exists some loop
$\xymatrix{
p\ar@(ur,r)^{z'_p}
}$ in $\A$ with $z'_p\pi = -n$.}$
\ei

Let 
$$r = \max(\{ m\} \cup \{ ||z_p|| \; \big{\lvert}\; p \in P^+\} \cup \{ ||z'_p|| \; \big{\lvert}\; p \in P^-\})
\quad\mbox{and}\quad
n' = 2(m^4n+r).$$
Since $P^+$ and $P^-$ are computable, also $r$ and $n'$ are computable.

We prove that (\ref{xix1}) holds for $i = 1$. We start by showing that
\beq
\label{afo7}
\mbox{if $u\xi_{n,1} \in G_{m+2r,1} \cap L$, then $u \in L$.}
\eeq
 
Let $x \in u\xi_{n,1}\theta\inv \cap L(\A)$, say $x = a_1\ldots a_s$ with $a_1,\ldots,a_s \in A$. Then there exists a path in $\A$ of the form
\beq
\label{afo5}
I \ni q_0 \xr{a_1} q_1 \xr{a_2} \ldots \xr{a_s} q_s \in T.
\eeq

There exists some $\ell \in [s]$ such that 
$$(a_1\ldots a_{\ell})\pi = u\xi_{n,1}\lambda \leq -m-2r.$$ 
Then there exist some $1 \leq j_0 <  \ldots < j_m < \ell$ such that
$(a_1\ldots a_{j_h})\pi = -r-h$ for $h = 0,\ldots,m$. We can assume that $j_0$ is largest possible. Similarly, there exist some $\ell < j'_m <  \ldots < j'_0 \leq s$ such that
$(a_1\ldots a_{j'_h})\pi = -r-h$ for $h = 0,\ldots,m$. We can assume that $j'_0$ is smallest possible.
Now there exist some $0 \leq b < b' \leq m$ such that $q_{j_b} = q_{j_{b'}}$, hence $p' = q_{j_b} \in P^-$. Similarly, there exist some $0 \leq c < c' \leq m$ such that $q_{j'_c} = q_{j'_{c'}}$, hence $p = q_{j_c} \in P^+$. 

Let 
$$y = a_1\ldots a_{j_b} z'_{p'} a_{j_b+1} \ldots a_{j'_c} z_{p} a_{j'_c+1} \ldots a_s \in L(\A).$$
We claim that $y\theta = x\theta$. Clearly,
$$y\pi = x\pi + z'_{p'}\pi + z_{p}\pi = x\pi -n+n = x\pi.$$
We have $y\rho \leq x\rho$ because the insertion of $z'_{p'}$ decreases $\pi$ and  $(a_1\ldots j_b)\pi \leq -r$ (and the insertion of $z_p$ cannot increase $\rho$ either). We actually have $y\rho = x\rho$ by maximality of $j_0$ and minimality of $j'_0$: this ensures that $x\rho$ is reached outside the interval $[j_0,j'_0]$. Finally, $y\lambda = x\lambda -n$ follows from the fact that $x\lambda$ is reached at $q_{\ell}$ and so we benefit from the insertion of $z'_{p'}$ which decreases $\pi$ by $n$. Thus $y\theta = x\theta = u$ and so $u \in (L(\A))\theta = L$. Therefore (\ref{afo7}) holds.

Now let $u \in G_{n',1} \subseteq G_{m+2r,1}$. If $u\xi_{n,1} \in L$, then $u \in L$ by (\ref{afo7}). 

Conversely, suppose that $u \in L$. Let $x \in u\theta\inv \cap L(\A)$, say $x = a_1\ldots a_s$ with $a_1,\ldots,a_s \in A$. Then there exists a path in $\A$ of the form
(\ref{afo5}).

Write
$$J = \{ j \in [s] \mid (a_1\ldots a_j)\pi = -m^2n\},\quad
J' = \{ j \in [s] \mid (a_1\ldots a_j)\pi = -2m^4n\}.$$
Since $n' > 2m^4n > m^2n$, both $J$ and $J'$ are nonempty. Consider the elements of $J \cup J'$ listed under the usual ordering. Since $u\pi \geq 0$, the first and last elements of the list belong to $J$. 
Now we remove from the list:
\bi
\item
all the elements of $J$ which are not adjacent to an element of $J'$;
\item
all the elements of $J'$ which are not adjacent to an element of $J$.
\ei
Thus our list takes the form
$$j_0 < j'_1 \leq j'_2 < j_3 \leq  j_4 < j'_5 \leq j'_6 < j_7 \leq \ldots  \leq j_{4k} < j'_{4k+1} \leq j'_{4k+2} < j_{4k+3},$$
with $j_{\ell} \in J$ and $j'_{\ell} \in J'$. Since $u\lambda < -2m^4n$, there exists some $t_0 \in \{ 0,\ldots,k\}$ such that $u\lambda = (a_1\ldots a_j)\pi$ with $j'_{4t_0+1} \leq j \leq j'_{4t_0+2}$. Indeed, any other possibility would imply the existence of another element of $J'$ that should not have been removed from the list!

There exist
$$j_{4t_0} = i_0 < i_1 < \ldots < i_{m^2} < j'_{4t_0+1} \leq j'_{4t_0+2} < i'_{m^2} < \ldots < i'_1 < i'_0 = j_{4t_0+3}$$
such that 
\beq
\label{afo6}
(a_1\ldots a_{i_{\ell}})\pi = (a_1\ldots a_{i'_{\ell}})\pi = -(m^2+\ell)n
\eeq
for $\ell = 0,\ldots,m^2$.  Hence there exist $0 \leq \ell_1 < \ell_2 \leq m^2$ such that $(q_{i_{\ell_1}}, q_{i'_{\ell_1}}) = (q_{i_{\ell_2}}, q_{i'_{\ell_2}})$. Now we remove from $x$ the factors $a_{i_{\ell_1}+1}\ldots a_{i_{\ell_2}}$ and $a_{i'_{\ell_2}+1}\ldots a_{i'_{\ell_1}}$ to get a word $x' \in L(\A)$. 

It follows from (\ref{afo6}) that $x'\pi = x\pi$. On the other hand, we removed before all those elements from the list $J \cup J'$ to ensure that removing these two factors does not alter $\rho$ because 
$(a_1\ldots a_h)\pi \leq -m^2n$ whenever $j_{4t_0} \leq h \leq j_{4t_0+3}$. 
What about $\lambda$? In what refers to the interval $[j_{4t_0}, j_{4t_0+3}]$, $\lambda$ would increase by $(\ell_2-\ell_1)n$. But we must control what happens in other intervals too, i.e. we must ensure that $\lambda$ will increase at least as much.

Thus let $t \in \{ 0,\ldots,k\}\setminus \{ t_0\}$. There exist
$$j_{4t} = i_0 < i_1 < \ldots < i_{m^2} < j'_{4t+1} \leq j'_{4t+2} < i'_{m^2} < \ldots < i'_1 < i'_0 = j_{4t+3}$$
such that 
$$(a_1\ldots a_{i_{\ell}})\pi = (a_1\ldots a_{i'_{\ell}})\pi = -m^2(1+\ell)n$$
for $\ell = 0,\ldots,m^2$. Similarly to the case of $t_0$, there exist $0 \leq \ell'_1 < \ell'_2 \leq m^2$ such that $(q_{i_{\ell'_1}}, q_{i'_{\ell'_1}}) = (q_{i_{\ell'_2}}, q_{i'_{\ell'_2}})$. Now we remove from $x$ the factors $a_{i_{\ell'_1}+1}\ldots a_{i_{\ell'_2}}$ and $a_{i'_{\ell'_2}+1}\ldots a_{i'_{\ell'_1}}$. Performing this operation for all $t$, we still get a word $x'' \in L(\A)$. 

For the same reasons invoked for $x'$, we have $x''\pi = x'\pi = x\pi$ and $x''\rho = x'\rho = x\rho$. In what refers to the interval $[j_{4t}, j_{4t+3}]$, $\lambda$ would increase by $(\ell'_2-\ell'_1)m^2n \geq (\ell_2-\ell_1)n$, hence $x''\lambda$ is reached in the interval $[j_{4t_0}, j_{4t_0+3}]$
and so $x''\lambda = x\lambda+(\ell_2-\ell_1)n$. Thus $u\xi_{n,1}^{\ell_2-\ell_1} = x''\theta \in L$. 

Now $u \in G_{n',1}$ implies 
$$u\xi_{n,1}^{\ell_2-\ell_1} \in G_{n'-(\ell_2-\ell_1)n,1} \subseteq G_{n'-m^2,1} \subseteq G_{m+2r,1}.$$
Applying (\ref{afo7}) $\ell_2-\ell_1-1$ times, we get $u\xi_{n,1} \in L$. Therefore (\ref{xix1}) holds also in this final case.

Assume now that $i = 3$. Since $L\beta \in \rat(F)$, we already know that there exist some computable $n,n' \geq 1$ such that 

\beq
\label{beta5}
u \in L\beta \iff u\xi_{n,1} \in L\beta
\eeq
holds for every $u \in G_{n',1}$.

Let $u \in G_{n',3}$. Then $u\beta \in G_{n',1}$. Since $\beta^2 = 1$, it follows from (\ref{beta4}) and (\ref{beta5}) that 
$$u \in L \iff u\beta \in L\beta \iff u\beta\xi_{n,1} \in L\beta \iff u\xi_{n,3}\beta \in L\beta \iff u\xi_{n,3} \in L,$$
hence (\ref{xix1}) holds for $i = 3$.

The case $i = 2$ is dealt with similarly to the case $i = 1$. This time we take
$$J = \{ j \in [s] \mid (a_1\ldots a_j)\pi = m^2n\},\quad
J' = \{ j \in [s] \mid (a_1\ldots a_j)\pi = 2m^4n\}.$$
The first element of the list built from $J \cup J'$ belongs to $J$, the last belongs to $J'$. We omit the details.
\qed

The following corollary will allow us to assume extra conditions on $n$ and $n'$.

\bc
\label{mul}
Let $L \in {\rm Rat}(F)$. For all $k > 1$ and $n'' \geq n'+kn$, we can replace the constants $n,n' \geq 1$ of Lemma \ref{xix} by $kn, n''$.
\ec

\proof
Let $i \in [3]$ and $u \in G_{n'',i}$. Since $n''-kn \geq n'$, we may apply (\ref{xix1}) $k$ times to get
$u\xi_{n,i}^{j-1} \in L \iff u\xi_{n_1,i}^j \in L$ for $j \in [k]$. By transitivity, we get 
$$u \in L \iff u\xi_{n,i}^{k} \in L \iff u\xi_{kn,i} \in L.$$
Thus (\ref{xix1}) holds for $kn,n''$.
\qed

\section{Deciding equality}

We prove in this section that inclusion and equality are decidable within $\rat(F)$.

Given $L \subseteq A^*$, let 
$$L^+ = L \cap \N\pi\inv\quad\mbox{and}\quad L^- = L \cap (-\N)\pi\inv.$$

Given $n,n' \geq 1$, we define a mapping $\eta_{n,n'}: \Z \to \Z$ as follows. Given $m \geq 0$, let
$$m\eta_{n,n'} = \left\{
\begin{array}{ll}
m&\mbox{ if }m < n'\\
\max\{ i \in \{ 0,n'-1\} \mid i \equiv m\,({\rm mod}\, n) \}&\mbox{ otherwise}
\end{array}
\right.
$$
and $(-m)\eta_{n,n'} = -(m\eta_{n,n'}).$

Let 
$$W_{n'} = \{ u \in F\mid u\lambda > -n',\, 0 \leq u\pi < n',\, u\rho-u\pi < n'\}.$$ 
We define also a mapping 
$\zeta_{n,n'}:F^+ \to W_{n'}$ as follows. Given $u \in F^+$, then
$$u\zeta_{n,n'}\lambda = (u\lambda)\eta_{n,n'}, \quad u\zeta_{n,n'}\pi = (u\pi)\eta_{n,n'},\quad u\zeta_{n,n'}\rho-u\zeta_{n,n'}\pi = (u\rho-u\pi)\eta_{n,n'}.$$
That is, $u\zeta_{n,n'}$ is obtained from $u$ by successively:
\bi
\item
applying $\xi_{n,1}$ while the value of $\lambda$ is $\leq -n'$;
\item
applying $\xi_{n,2}$ while the value of $\pi$ is $\geq n'$;
\item
applying $\xi_{n,3}$ while the value of $\rho-\pi$ is $\geq n'$.
\ei
It is immediate that $\zeta_{n,n'}^2 = \zeta_{n,n'}$.

\bt
\label{form}
Let $L \in {\rm Rat}(F)$. Then there exist computable $n' \geq n \geq 1$ such that:
\bi
\item[(i)] $L^+ = L^+\zeta_{n,n'}\zeta_{n,n'}\inv = (L \cap W_{n'})\zeta_{n,n'}\inv$;
\item[(ii)] $L^- = ((L^-)\inv\zeta_{n,n'}\zeta_{n,n'}\inv)\inv = ((L\inv \cap W_{n'})\zeta_{n,n'}\inv)\inv$.
\ei
\et

\proof
First we note that $L\inv \in \rat(F)$: if $\A = (Q,I,T,E)$ is a finite $A$-automaton such that $(L(\A))\theta = L$, then $L\inv = (L(\A'))\theta$ for the automaton $\A'$ obtained from $\A$ by replacing each edge $p \xr{a^{\varepsilon}} q$ by $q \xr{a^{-\varepsilon}} p$ and exchanging the initial with the terminal vertices.

Let $n'_1 \geq n_1 \geq 1$ (respectively $n'_2 \geq n_2 \geq 1$) be the (computable) constants given in Lemma \ref{xix} for $L$
(respectively $L\inv$). Let $n = n_1n_2$ and let $n' = \max\{n'_1,n'_2\} + n$. By Corollary \ref{mul}, we may replace both $n_1,n'_1$ and $n_2,n'_2$ by $n,n'$.

Write $\eta = \eta_{n,n'}$ and $\zeta = \zeta_{n,n'}$.

Let $u \in L^+\zeta\zeta\inv$. Then $u\zeta = v\zeta$ for some $v \in L^+$. 
There exist $i,j,k \geq 0$ such that $v\zeta = v\xi_{n,1}^i\xi_{n,2}^j\xi_{n,3}^k$. Suppose that $i > 0$. Then $(v\lambda)\eta = v\zeta\lambda > v\lambda$ and so $-n' < (v\lambda)\eta \leq -(n'-n)$. Thus $v \in G_{n',1}$. Similarly, $j > 0$ implies $v \in G_{n',2}$, and $k > 0$ implies $v \in G_{n',3}$. Therefore we may apply
 Lemma \ref{xix} and get $v\zeta \in L$. 
 
 On the other hand, $v\zeta = u\zeta = u\xi_{n,1}^{i'}\xi_{n,2}^{j'}\xi_{n,3}^{k'}$ for some $i',j',k' \geq 0$. The same argument allows us to use Lemma \ref{xix} and get $u \in L$. Thus $u \in L \cap F^+ = L^+$.
Therefore $L^+\zeta\zeta\inv \subseteq L^+$. Since the opposite inclusion holds trivially, we obtain $L^+ = L^+\zeta\zeta\inv$.

To complete the proof of (i), it suffices to show that
$$L^+\zeta = L \cap W_{n'}.$$

Let $u \in L^+$. Clearly, $u\zeta \in W_{n'}$. On the other hand, 
$u\zeta^2 = u\zeta \in L^+\zeta$, hence $u\zeta \in L^+\zeta\zeta\inv = L^+$ and so $L^+\zeta \subseteq L^+ \cap W_{n'}  \subseteq L \cap W_{n'}$.

Conversely, let $u \in L \cap W_{n'}$. Then $u \in L^+$ and so $u = u\zeta \in L^+\zeta$. Thus $L^+\zeta = L \cap W_{n'}$ and (i) holds.

Since $(L^-)\inv = (L\inv)^+$, we get (ii) by applying (i) to $L\inv \in \rat(F)$.
\qed

We can now prove that inclusion and equality are decidable in $\rat(F)$.

\bt
\label{dec}
Let $K,L \in {\rm Rat}(F)$. Then it is decidable whether or not:
\bi
\item[(i)] $K \subseteq L$;
\item[(ii)] $K = L$.
\ei
\et

\proof
(i) In view of Lemma \ref{xix} and Corollary \ref{mul}, there exist computable constants $n' \geq n \geq 1$ such that (\ref{xix1}) holds for both $K$ and $L$.
Let $\zeta = \zeta_{n,n'}$. We claim that
\beq
\label{dec1}
K^+ \subseteq L^+\quad\mbox{if and only if}\quad K \cap W_{n'} \subseteq L \cap W_{n'}.
\eeq

The direct implication follows from $W_{n'} \subseteq F^+$. The converse follows from Theorem \ref{form}(i). Thus  (\ref{dec1}) holds.

Since $K\cap W_{n'}$ and $L\cap W_{n'}$ are finite and effectively computable in view of Corollary \ref{mpr}, we can decide whether or not $K^+ \subseteq L^+$.

Now
$$K \subseteq L \iff (K^+ \subseteq L^+ \,\wedge K^- \subseteq L^-) \iff (K^+ \subseteq L^+ \,\wedge (K\inv)^+ \subseteq (L\inv)^+),$$
 hence we can decide whether or not $K \subseteq L$.

(ii) Follows from (i).
\qed

\bc
\label{decmon}
Let $L \in {\rm Rat}(F)$. Then it is decidable whether or not $L$ is a submonoid of $F$.
\ec

\proof
Indeed, $L$ is a submonoid of $F$ if and only if $L^* = L$, and this equality is decidable by Theorem \ref{dec}.
\qed

Theorem \ref{form} also yields the following corollary:

\bc
\label{mai}
Let $L \in {\rm Rat}(F)$. Then there exist computable $n' \geq n \geq 1$ and $W,W' \subseteq W_{n'}$ such that $L =  W\zeta_{n,n'}\inv \cup (W'\zeta_{n,n'}\inv)\inv$. 
\ec

Could this provide a characterization of $\rat(F)$?
The following example settles the question in the negative, showing also that $\rat(F)$ is not closed under complementation. 

\be
\label{ncuc}
Let $L = Faa\inv$ and $W = \{ 1,a,a\inv a, a\inv a^2\} \subseteq W_2$. Then $L \in {\rm Rat}(F)$ but $F\setminus L = W\zeta_{1,2} \notin {\rm Rat}(F)$.
\ee

Indeed, $L \in \rat(F)$ obviously. Suppose that $F \setminus L \in \rat(F)$. Then $F\setminus L = (L(\A))\theta$ for some finite $A$-automaton $\A = (Q,I,T,E)$. Let $m = |Q|$ and $u = a^{-m}a^m \in A^+$. Since $u\theta \in F \setminus L$, there is some path $I \ni q_0 \xr{v} t \in T$ in $\A$ with $v\theta = u\theta$. 

Since $v\theta\lambda = -m$, there exists some factorization $v = v_1\ldots v_mv'$ such that 
$v_i\theta\pi = -1$ for $j \in [m]$. And we have a path
$$I \ni q_0 \xr{v_1} q_1 \xr{v_2} \ldots \xr{v_m} q_m \xr{v'} t \in T$$
in $\A$. Since $m = |Q|$, there exist some $0 \leq i < j \leq m$ such that $q_i = q_j$. Hence $(v_{i+1}\ldots v_j)\theta\pi = i-j < 0$. Let $w = v_1\ldots v_i(v_{i+1}\ldots v_j)^2v_{j+1}\ldots v_mv' \in L(\A)$. Since $w\pi = i-j < 0$, we get $w\theta \in L$, contradicting $w\theta \in (L(\A))\theta = F\setminus L$. Therefore $F \setminus L \notin \rat(F)$ as claimed.

The elements of $F\setminus L$ are those $u \in F$ which have a Munn tree of the form
$$\xymatrix{
\cdot \ar[rrr]^{a^i} &&& \circ \ar[rr]^{a^j} && \bullet
}$$
for some $i,j \geq 0$. It is easy to check that
$$u\zeta_{1,2} = \left\{
\begin{array}{ll}
1&\mbox{ if }i = j = 0\\
a&\mbox{ if }i = 0 < j\\
a\inv a&\mbox{ if }i > 0 = j\\
a\inv a^2&\mbox{ if }i,j > 0
\end{array}
\right.$$
and $F\setminus L = W\zeta_{1,2}$.

\section{Recognizable subsets}

We prove in this section that it is decidable whether or not a given rational subset of $M$ is recognizable.

Given a subset $L$ of a monoid $M$, we define the {\em syntactic congruence} of $L$ as the congruence $\sim_L$ on $M$ defined by
$$u\sim_L v\quad\mbox{if}\quad \forall x,y \in M\,(xuy \in L \iff xvy \in L).$$
We say that $L$ is {\em recognizable} if $\sim_L$ has finite index (i.e. if $M/\!\sim_L$ is a finite monoid). Equivalently, $L$ is recognizable if there exists some homomorphism $\p$ from $M$ onto a finite monoid $N$ such that $L = L\p\p\inv$.
We denote by $\rec(M)$ the set of all recognizable subsets of $M$. If $M$ is finitely generated, then $\rec(M) \subseteq \rat(M)$. See \cite{Ber} for more details on recognizable languages.

The next example shows that $\rec(F) \subset \rat(F)$:

\be
Let $L = a^* \subset F$. Then $L \in {\rm Rat}(F) \setminus {\rm Rec}(F)$.
\ee

Obviously, $L \in \rat(F)$, so it suffices to show that $a^m \sim_L a^n$ implies $m = n$ for all $m,n \geq 0$. Indeed, if $m< n$, then $a^ma^{-n}a^n = a^{-(n-m)}a^n \notin L$ but $a^na^{-n}a^n = a^n \in L$.

\bt
\label{crec}
Let $L \in {\rm Rat}(F)$. Let $n,n' \geq 1$ and $W,W' \subseteq W_{n'}$ be such that $L^+ =  W\zeta_{n,n'}\inv$ and  $L^- = (W'\zeta_{n,n'}\inv)\inv$.
Then the following conditions are equivalent:
\bi
\item[(i)] $L \in {\rm Rec}(F)$;
\item[(ii)] $a^{3n'+n} \sim_L a^{3n'}$.
\ei
\et

\proof
(i) $\Rw$ (ii).
Write $\zeta = \zeta_{n,n'}$ and $\tau = \;\sim_L$. 
Since $L \in \rec(F)$, there exist $m,p \geq 1$ such that $a^{m} \,\tau\, a^{m+p}$. 
We show that
\beq
\label{crec1}
\mbox{if $||w|| \geq m$, then there exists $k_0 \in \N$ such that $w \,\tau\, a^{w\pi+kp}$ for every $k \geq k_0$.}
\eeq

Indeed, take some $k_0$ such that $w\pi + k_0p \geq m$. Let $k \geq k_0$ and write $u = a^{w\pi+kp}$. 
Since $||w|| \geq m$, then $w = w_1a^mw_2$ for some $w_1,w_2 \in F$. Let $v = w_1a^{m+kp}w_2$. Since $a^{m} \,\tau\, a^{m+p}$, we get $a^{m} \,\tau\, a^{m+kp}$ and consequently $w \,\tau\, v$. Now $v = vv\inv a^{v\pi} = vv\inv u$.

By \cite[Lemma 4.6]{Lee}, the minimal ideal of any finite monogenic inverse monoid is a cyclic group. Let $K$ denote 
the minimal ideal of $A^*/\tau$. In view of the relation $a^m \,\tau\, a^{m+p}$, we have necessarily $a^m\tau \in K$.
Since $||u||, ||v|| \geq m$, we have $u\tau, v\tau \in K$. Hence
$w\tau = v\tau = (vv\inv u)\tau = u\tau$ and (\ref{crec1}) holds.

We show next that if $w \in F^+$ is such that $||w|| \geq 3n'$, then there exists some $\oo{w} \in F^+$ such that
$$||\oo{w}|| \geq m,\quad \oo{w}\zeta = w\zeta \quad\mbox{and}\quad
w\pi \equiv \oo{w}\pi \,({\rm mod}\, n).$$

Since $||w|| \geq 3n'$, then $w \in G_{n',i}$ for some $i \in [3]$. It is straightforward to check that all the conditions are satisfied by $\oo{w} = w\xi_{n,i}^m$.

By taking inverses, we see also that, if $w \in F^-$ is such that $||w|| \geq 3n'$, then there exists some $\oo{w} \in F^-$ such that
$$||\oo{w}|| \geq m,\quad \oo{w}\inv\zeta = w\inv\zeta \quad\mbox{and}\quad
w\pi \equiv \oo{w}\pi \,({\rm mod}\, n).$$

The next step is to show that
\beq
\label{crec2}
\mbox{if $||u||, ||v|| \geq m$ and $u\pi \equiv v\pi \,({\rm mod}\, n)$, then $u \in L$ if and only if $v\in L$.}
\eeq

We may assume that $u\pi > v\pi$, hence $u\pi = v\pi + rn$ for some $r \geq 1$. By (\ref{crec1}), there exists some $k \in \N$ such that $v\pi+kp \geq n'$, $u \,\tau\, a^{u\pi+kp}$ and $v \,\tau\, a^{v\pi+kp}$. Now $a^{u\pi+kp}\zeta = a^{v\pi+kp+rn}\zeta = a^{v\pi+kp}\zeta$, hence 
$$u \in L \iff a^{u\pi+kp} \in L \iff a^{u\pi+kp}\zeta \in W \iff a^{v\pi+kp}\zeta \in W \iff a^{v\pi+kp} \in L \iff v\in L$$
and (\ref{crec2}) holds.

Now let $x,y \in F$. Since 
$||xa^{3n'+n}y||, ||xa^{3n'}y|| \geq 3n'$, we have $||\oo{xa^{3n'+n}y}||, ||\oo{xa^{3n'}y}|| \geq m$ and it follows from (\ref{crec2}) that 
\beq
\label{crec5}
\oo{xa^{3n'+n}y}\in L\mbox{ if and only if } \oo{xa^{3n'}y} \in L.
\eeq
On the other hand, we have
$$xa^{3n'+n}y \in L^+ \iff (xa^{3n'+n}y)\zeta \in W \iff \oo{xa^{3n'+n}y}\zeta \in W \iff \oo{xa^{3n'+n}y} \in L^+$$
and
$$xa^{3n'+n}y \in L^- \iff (xa^{3n'+n}y)\inv\zeta \in W' \iff \oo{xa^{3n'+n}y}\inv\zeta \in W' \iff \oo{xa^{3n'+n}y} \in L^-,$$
hence $xa^{3n'+n}y \in L$ if and only if $\oo{xa^{3n'+n}y} \in L$. Similarly, $xa^{3n'}y \in L$ if and only if $\oo{xa^{3n'}y} \in L$.

Together with (\ref{crec5}), this implies that
$$xa^{3n'+n}y \in L \iff \oo{xa^{3n'+n}y}\in L \iff \oo{xa^{3n'}y} \in L \iff xa^{3n'}y \in L.$$
Therefore $a^{3n'+n} \sim_L a^{3n'}$.

(ii) $\Rw$ (i). It is immediate that, for every $w \in F^+$, there exists some $w' \in F^+$ such that $||w'|| \leq 3(3n'+n)$ and $w' \tau w$. By taking inverses, this also holds for $w \in F^-$, therefore $F/\tau$ is finite and so $L \in \rec(M)$.
\qed

We can now prove the announced decidability result:

\bt
\label{dere}
Given $L \in {\rm Rat}(F)$, it is decidable wether or not $L \in {\rm Rec}(F)$.
\et

\proof
By Corollary \ref{mai}, we can compute $n,n' \geq 1$ and $W,W' \subseteq W_{n'}$ such that $L =  W\zeta_{n,n'}\inv \cup (W'\zeta_{n,n'}\inv)\inv$. By Theorem \ref{crec}, it suffices to show that it is decidable whether or not $a^{3n'+n} \sim_L a^{3n'}$.

Let $M$ be the inverse monoid defined by the inverse monoid presentation $\la a \mid a^{3n'+n} = a^{3n'}\ra$ (i.e., $M$ is the quotient of $F$ by the congruence generated by $(a^{3n'+n},a^{3n'})$). By the argument used to prove the second implication of Theorem \ref{crec}, $M$ is finite. Let $\p:F \to M$ and $\pi:F/\!\sim_L$ be the canonical homomorphisms. We show that
\beq
\label{dere1}
a^{3n'+n} \sim_L a^{3n'} \mbox{ if and only if }L = L\p\p\inv.
\eeq

Suppose first that $a^{3n'+n} \sim_L a^{3n'}$. Then $a^{3n'+n}\pi = a^{3n'}\pi$ and so there exists a homomorphism $\psi:M \to F/\!\sim_L$ such that $\p\psi = \pi$. 

Let $u \in L\p\p\inv$. Then $u\p = v\p$ for some $v \in L$ and so $u\pi = u\p\psi = v\p\psi = v\pi$. Thus $u \sim_L v$. Since $v \in L$, this yields $u \in L$. Therefore $L = L\p\p\inv$.

Conversely, $L = L\p\p\inv$ implies $L \in \rec(F)$ because $M$ is finite and so $a^{3n'+n} \sim_L a^{3n'}$ in view of Theorem \ref{crec}. Therefore (\ref{dere1}) holds.

Thus it suffices to show that the equality $L = L\p\p\inv$ is decidable. But $L\p$ is a computable subset of the finite monoid $M$ and so $L\p\p\inv$ is a computable recognizable subset of $F$. Since $F$ is finitely generated, we get $L\p\p\inv \in \rat(F)$. Therefore we can decide $L = L\p\p\inv$ by Theorem \ref{dec} and we are done.
\qed

\section{Rational submonoids}

The main goal of this section is to prove a version of the Anisimov and Seifert Theorem valid for the monogenic free inverse monoid $F$.

We start by proving a couple of lemmas.

\bl
\label{pone}
Let $M$ be a rational submonoid of $F$ such that $M\pi$ contains both positive and negative integers. Then $M\pi = n\Z$ for some computable $n \in \N \setminus \{ 0\}$.
\el

\proof
The subgroup of $\Z$ generated by $M\pi$ is necessarily cyclic, say $n\Z$, for some $n \geq 1$. Rational subgroups of  free groups are finitely generated by Anisimov and Seifert's theorem \cite[Theorem 4.1]{BS1} and a finite generating set can always be computed. Thus $n$ is computable.

Now there exist $m_1, \ldots, m_k \in M\pi\setminus \{ 0\}$ and $x_1,\ldots,x_k \in \Z$ such that $n = m_1x_1 + \ldots + m_kx_k$. Without loss of generality, we may assume that $m_1 > 0$ and $m_k < 0$ because some $x_i$ can be chosen to be zero. Now assume that the number $j$ of negative $x_i$ is minimum. Suppose that $j > 0$. Take some $x_i < 0$.

Suppose that $m_i < 0$. By replacing $x_1$ by $x_1+m_ix_i$ and $x_i$ by $x_i-m_1x_i$, we get an alternative decomposition of $n$ which contradicts the minimality of $j$ since $x_1+m_ix_i > x_1$ and $x_i-m_1x_i \geq 0 > x_i$. 

Hence $m_i > 0$. Now we replace instead $x_i$ by $x_i+m_kx_i$ and $x_k$ by $x_k-m_ix_i$, yielding an alternative decomposition of $n$ which contradicts the minimality of $j$ since $x_i+m_kx_i \geq 0 > x_i$ and $x_k-m_ix_i > x_k$. 
Thus $j = 0$ and so $n \in M\pi$.

It is easy to see that $M\inv$ is itself a rational monoid. In what concerns rationality, this follows from
$$(X\cup Y)\inv = X\inv \cup Y\inv, \quad (XY)\inv = Y\inv X\inv, \quad (X^*)\inv  = (X\inv)^*.$$
Clearly, the subgroup of $\Z$ generated by $M\inv\pi$ is still $n\Z$, so it follows from what we have already proved that $n \in M\inv\pi$. Thus $-n \in M\pi$ and so $M\pi = n\Z$.
\qed

\bl
\label{gon}
Let $L \in {\rm Rat}(F)$. Then there exists some computable $K \in \N$ such that, for all $s \geq 1$ and $u \in L^+$ such that $u\rho-u\pi \geq K$, there exists some $v \in L$ such that 
\beq
\label{cab1}
v\lambda -v\pi \geq u\lambda-u\pi,\quad -s-K < v\pi \leq -s\quad \mbox{and}\quad v\rho = u\rho.
\eeq 
\el

\proof
Let $\A = (Q,I,T,E)$ be a finite $A$-automaton such that $(L(\A))\theta = L$. 
Let $K = |Q|$.
Fix $x \in u\theta\inv \cap L$. Then there exists a path in $\A$ of the form
\beq
\label{gon1}
I \ni q_0 \xr{x'} p_0 \xr{x_1} p_1 \xr{x_2} \ldots \xr{x_K} p_{K} \xr{x''} t \in T
\eeq
with $x = x'x_1\ldots x_{K}x''$, $x'\pi = x\pi+K$ and $x_r\pi = -1$ for $r = 1,\ldots, K$. 
We may assume that $x',x_1,\ldots,x_{K}$ are successively chosen as having maximum length. This amounts to say that the enhanced occurrences of each $p_r$ in (\ref{gon1}) are the last occurrences of each such vertex in the path.

There exist $0 \leq r_1 < r_2 \leq K$ such that $p_{r_1} = p_{r_2}$. Let
$y = x_1\ldots x_{r_1}$, $y' = x_{r_1+1}\ldots x_{r_2}$ and $y'' = x_{r_2+1}\ldots x_K x''$. Then $x = x'yy'y''$. Since $r_2-r_1 \in [K]$, there exists some $n \in \N$ such that 
$$-s-K < u\pi - (n-1)(r_2-r_1) \leq -s.$$
Then $x'y(y')^ny'' \in L(\A)$. We define $v = (x'y(y')^ny'')\theta \in L$. Note that
$$v\pi = u\pi + (y')^{n-1}\pi = u\pi -(n-1)(r_2-r_1) \in \{ -s-K+1,\ldots,-s\}.$$ 

The following picture depicts the general situation in terms of Munn trees, where the dotted lines are potential lines: 

$$
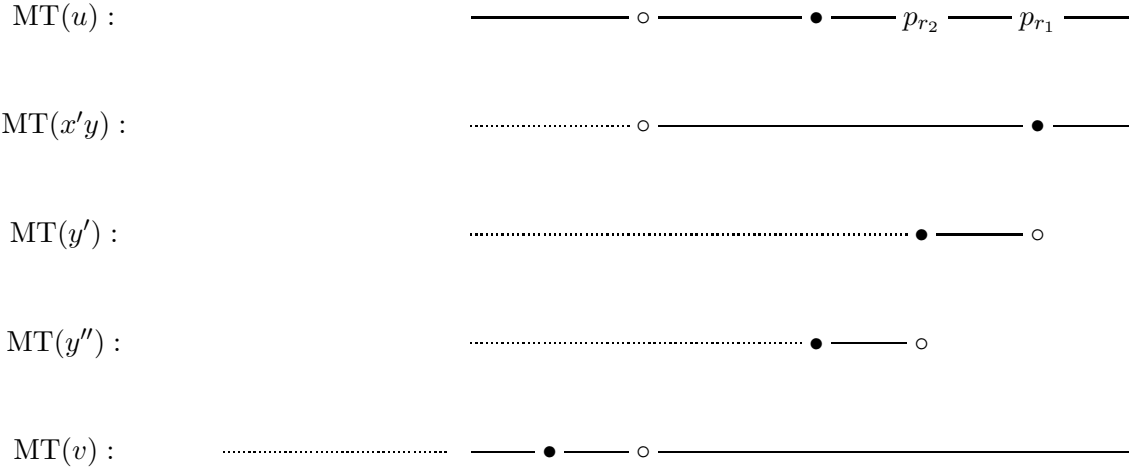
\xymatrix{
\mt(u):&&&& \ \ar@{-}[rr] && \circ \ar@{-}[rr] && \bullet \ar@{-}[r] & p_{r_2} \ar@{-}[r] & p_{r_1} \ar@{-}[r] &\\
\mt(x'y):&&&& \ \ar@{.}[rr] && \circ \ar@{-}[rrrr] &&&& \bullet \ar@{-}[r] &\\
\mt(y'):&&&& \ \ar@{.}[rrrrr] &&&&& \bullet \ar@{-}[r] & \circ &\\
\mt(y''):&&&&\ \ar@{.}[rrrr] &&&& \bullet \ar@{-}[r] & \circ &&\\
\mt(v):&\ \ar@{.}[rrr] &&& \ \ar@{-}[r] & \bullet \ar@{-}[r] & \circ \ar@{-}[rrrrr] &&&&&
}$$
Note that the shapes of $\mt(x'y)$, $\mt(y')$ and $\mt(y'')$ follow from the enhanced occurrences of each $p_r$ in (\ref{gon1}) being the last occurrences of each such vertex in the path. And $u\rho = x'\rho = v\rho$ follows from $x',x_1,\ldots,x_{K}$ having maximum length.

Finally, note that $u\lambda-u\pi = (x'yy'y'')\lambda - (x'yy'y'')\pi$. Each insertion of $y'$ decreases $\pi$ by $y'\pi$ and $\lambda$ by at most $y'\pi$. Thus $v\lambda -v\pi \geq u\lambda-u\pi$ and we are done.
 \qed

We can now prove the main result of this section. We will consider some cases and subcases.

\bt
\label{fg}
Let $M$ be a rational submonoid of $F$. Then $M$ is finitely generated and we can compute a finite generating set for $M$.
\et

\proof
Let $\A = (Q,I,T,E)$ be a finite $A$-automaton such that $(L(\A))\theta = L$. 
After introducing $\beta$, we remarked that there exists an automaton $(Q,T,I,E')$ recognizing $M\beta$. Thus $M\beta$ is also a rational submonoid of $F$. By the proof of Lemma \ref{gon}, the same constant $K = |Q|$ works for both $M$ and $M\beta$ if we need to apply the lemma. For technical reasons, it is sometimes convenient for us to deal simultaneously with $M$ and $M\beta$. Note that $\beta\pi = \pi$.

Let $n' \geq n \geq 1$ be the constants provided by Lemma \ref{xix} for $M$. By Corolary \ref{mul}, we may assume that $n,n'$ hold for $M\beta$ as well.

We consider three cases. It is easy to know which case we are in since $M\pi$ is clearly computable. 

\medskip

\noindent
\underline{Case 1}: $M\pi = \{0\}$.

\medskip

By Lemma \ref{gon}, we have $u\rho < K$ for every $u \in M$. Since $\beta\pi = \pi$, we may also apply Lemma \ref{gon} to $M\beta$ and get $u\beta\rho < K$. By (\ref{beta1}), we get $u\beta\rho = -u\lambda$, hence $u\lambda > -K$. Thus
$M$ is finite and consequently finitely generated. Computability follows from Corollary \ref{mpr}.

\medskip

\noindent
\underline{Case 2}: $M\pi \neq \{0\}$, and either $M = M^-$ or $M = M^+$.

\medskip

Replacing $M$ by $M\inv$ if necessary, we may assume that $M = M^+$. 
We consider
$$X = \{ x \in M \;\big{\lvert}\; ||x|| \leq (n+1)n' + n^2+ 2K\}.$$
This is a finite computable set by Corollary \ref{mpr}. Hence it suffices to show that $M \subseteq X^*$. Since $X^*\beta = (X\beta)^*$, this is equivalent to have $M\beta \subseteq (X\beta)^*$.
We prove both inclusions simultaneously by induction on the norm of $v \in M \cup M\beta$.

Assume then that $k > (n+1)n' + n^2+2K$ and:
\bi
\item
if $w \in M$ and $||w|| < k$, then $w \in X^*$;
\item
if $w \in M\beta$ and $||w|| < k$, then $w \in (X\beta)^*$.
\ei
We have to show that:
\bi
\item[(a)] if $v \in M$ and $||v|| = k$, then $v \in X^*$;
\item[(b)] if $v \in M\beta$ and $||v|| = k$, then $v \in (X\beta)^*$.
\ei
Since $\beta$ is an involution of $F$, then we could interchange (a) and (b) by replacing $M$ by $M\beta$. It is therefore
enough to consider $v \in M$ with $||v|| = k$.

Adapting the argument of Case 1, we have 
\beq
\label{ett1}
z\rho-z\pi < K\quad\mbox{and}\quad
v\lambda > -K
\eeq
for every $z \in M$.
Hence $v\pi > (n+1)n' +n^2$. 
Let 
$$r = \max\{ i \geq 1 \mid v\pi-in \geq n'\}.$$
By (\ref{xix1}), we have $w = v\xi_{n,2}^r \in M$. Also $n' \leq w\pi < n'+n$.
Let $z = v\xi_{n,2}^{w\pi}$. Then
$$v\pi -n(w\pi) > (n+1)n' +n^2 -n(n'+n) = n'$$
and it follows easily from (\ref{xix1}) that $z \in M$. Since $||w||,||z|| < k$, we have $w,z \in X^*$. Note that $w\lambda = z\lambda = v\lambda$ and $w\rho-w\pi = z\rho-z\pi = v\rho -v\pi$. On the other hand, 
$(w^{n-1}zw)\pi = n(w\pi) + z\pi = v\pi$. Since $w,z \in M^+$, we get $v = w^{n-1}zw \in X^*$ as desired.

\medskip

\noindent
\underline{Case 3}: $M\pi$ contains both positive and negative integers.

\medskip

By Lemma \ref{pone}, we have $M\pi = p\Z$ for some computable $p \in \N \setminus \{ 0\}$. Let $u,u' \in M$ be such that $u\pi = p$ and $u'\pi = -p$. Then $uu' \in E(M)$ and it follows easily that $u'uu' = (uu'u)\inv$. Since $(uu'u)\pi = p$, we can replace $u$ by $uu'u$ and assume that $u,u\inv \in M$.

By Corollary \ref{mul}, we may assume that $p \mid n$ and $n' \geq ||u||+n, K$.
We consider
$$X = \{ x \in M \;\big{\lvert}\; ||x|| \leq 18n'\}.$$
This is a finite computable set by Corollary \ref{mpr}. Hence it suffices to show that $M \subseteq X^*$. Out of symmetry, it is enough to show that $M^+ \subseteq X^*$. 
Since $X^*\beta = (X\beta)^*$, this is equivalent to have $M^+\beta \subseteq (X\beta)^*$.
We prove both inclusions simultaneously by induction on the norm of $v \in M^+ \cup M^+\beta$.

Assume then that $k > 18n'$ and:
\bi
\item
if $w \in M^+$ and $||w|| < k$, then $w \in X^*$;
\item
if $w \in M^+\beta$ and $||w|| < k$, then $w \in (X\beta)^*$.
\ei
By a similar reduction to the one we performed in Case 2, it is enough to consider $v \in M^+$ with $||v|| = k$.
We split the discussion into three subcases, since
$$||v|| = k > 18n' \Rw v \in G_{8n',1} \cup G_{2n',2} \cup G_{8n',3}.$$

\medskip

\noindent
\underline{Subcase 1}: $v \in G_{2n',2}$.

\medskip

Let $r = \min\{ i \geq 1 \mid v\pi-in \leq n'\}$. 
By (\ref{xix}), we have $w = v\xi_{n,2}^r \in M^+$. Since $||w|| < k$, we have $w \in X^*$. Since $M\pi = p\Z$, we have $v\pi -2(w\pi) = jp$ for some $j \geq 0$. We claim that $wu^jw = v$.

Indeed, $(wu^jw)\pi = 2(w\pi)+jp = v\pi$. On the other hand, 
$$-u\lambda \leq ||u|| \leq n'-n \leq w\pi$$
yields $(wu^jw)\lambda = w\lambda = v\lambda$. Similarly, 
$$u\rho-u\pi \leq ||u|| \leq n'-n \leq w\pi$$
yields $(wu^jw)\rho-(wu^jw)\pi = w\rho-w\pi = v\rho-v\pi$. Thus $v = wu^jw \in X^*$.

\medskip

\noindent
\underline{Subcase 2}: $v \in G_{8n',3}$.

\medskip

We may assume that $v\pi < 2n'$, otherwise we are covered by Subcase 1.
Let 
$$r = \max\{ i \geq 1 \mid v\rho-v\pi-in \geq n'\}.$$
By (\ref{xix}), we have $w = v\xi_{n,3}^r \in M^+$. Also $n' \leq w\rho-w\pi = v\rho-rn-v\pi < n'+n$.
Since $||w|| < k$, we have $w \in X^*$. 

Since $v\rho > 2n'$, there exists a path in $\A$ of the form
$$I \ni q_0 \xr{x_1} q_1 \xr{x_2} \ldots \xr{x_{2n'}} q_{2n'} \xr{x'} t \in T$$
with $v = (x_1\ldots x_{2n'}x')\theta$ and $x_j\pi = x_j\rho = 1$ for each $j \in [2n']$.
Then there exists some path $q_{2n'} \xr{x''} t$ with $|x''| < K \leq n'$. Let $v' = (x_1\ldots x_{2n'}x'')\theta$. It it is easy to check that
\beq
\label{rain1}
v' \in M^+, \quad v'\lambda \geq v\lambda\quad \mbox{and}\quad
n' \leq v'\pi \leq v'\rho < 2n'+K.
\eeq

On the other hand, it follows from (\ref{xix1}) that there is some $s \geq 0$ such that $w' = w\xi_{n,1}^s \in M^+$ and $w'\lambda > -n'$. By Lemma \ref{gon}, there exists some $w'' \in M$ such that 
\beq
\label{rain2}
w''\lambda -w''\pi \geq w'\lambda-w'\pi,\quad -n'-K < w''\pi \leq -n'\quad \mbox{and}\quad w''\rho = w'\rho.
\eeq 
We have $||v'|| = v'\rho-v'\lambda < 2n'+K-v\lambda < ||v|| = k$ and so $v' \in X^*$.
On the other hand,
$$\begin{array}{lll}
||w''||&=&w''\rho-w''\lambda \leq w'\rho -w'\lambda +w'\pi-w''\pi < w\rho + n' +w\pi +n' + K\\ &&\\
&<&2(w\pi) +5n' = 2(v\pi) + 5n' < 9n',
\end{array}$$
hence $w'' \in X$. Thus it suffices to show that 
\beq
\label{rain3}
v = wv'u^jwu^{-\ell}w'' \mbox{ for some }j,\ell \geq 0.
\eeq 

Let $j = \frac{v\rho-w\pi-v'\pi-w\rho}{p}$. Since $n$ divides $v\rho-w\rho$ and $M\pi = p\Z$, we have $j \in \Z$. On the other hand, we have
that
$$v\rho-w\pi-v'\pi-w\rho > v\rho-v\pi -2n'-K -n'-n-w\pi \geq 8n'-7n',$$
hence $jp > n'$ and $j > 0$. We claim that
\beq
\label{rain3}
v\lambda = (wv'u^jw)\lambda \quad\mbox{and}\quad
v\rho = (wv'u^jw)\rho.
\eeq 

Indeed, $w\lambda = v\lambda$ and $v',u,w \in M^+$. Since $v'\lambda \geq w\lambda$ and $u\lambda \geq -n' \geq -v'\pi$, we get $v\lambda = (wv'u^jw)\lambda$.

On the other hand,
$$(wv'u^jw)\rho \geq (wv'u^j)\pi + w\rho = w\pi +v'\pi +jp +w\rho = v\rho.$$
Now 
$$v'\rho-v'\pi \leq n'+K \leq 2n' < jp+w\rho \leq (u^jw)\rho$$
and $u\rho \leq n' \leq w\rho$.
Since $v',u,w \in M^+$, we must have $v\rho = (wv'u^jw)\rho$. Therefore (\ref{rain3}) holds.

Take
$$\ell = j+ \frac{w\pi+v'\pi +w''\pi}{p}.$$
Since $M\pi = p\Z$, we have $\ell \in \Z$. On the other hand, we have
$$\ell > \frac{n'+n'-(n'+K)}{p} \geq 0$$
and $v\pi = (wv'u^jwu^{-\ell}w'')\pi$.

By (\ref{rain3}), we have $(wv'u^jwu^{-\ell}w'')\rho \geq v\rho$. Since $u\inv \in M^-$ and $u\inv\rho \leq n' \leq w\rho$, we get 
$(wv'u^jwu^{-\ell})\rho = v\rho$. Finally, we note that
$$\begin{array}{lll}
(wv'u^jwu^{-\ell})\pi + w''\rho&=&2(w\pi)+v'\pi +(j-\ell)p+w'\rho =
2(v\pi)+v'\pi -(v\pi+v'\pi+w''\pi)+w\rho\\ &&\\
&=&v\pi-w''\pi+w\rho < v\pi+n'+K+n'+n+v\pi < 8n' < v\rho,
\end{array}$$
hence $(wv'u^jwu^{-\ell}w'')\rho = v\rho$.

It is straightforward to check that (\ref{rain3}) yields $v\lambda = (wv'u^jwu^{-\ell})\lambda$. Now our preceding computation yields
$$(wv'u^jwu^{-\ell})\pi + w''\lambda = w\pi -w''\pi + w''\lambda \geq w\pi -w'\pi + w'\lambda \geq w\pi -w\pi + w\lambda =
w\lambda = v\lambda,$$
thus $(wv'u^jwu^{-\ell}w'')\lambda = v\lambda$ and consequently $v = wv'u^jwu^{-\ell}w'' \in X^*$.

\medskip

\noindent
\underline{Subcase 3}: $v \in G_{8n',1}$.

\medskip

Then $v\beta \in G_{8n',3}$. In view of the induction hypothesis and Subcase 2 (with respect to $M\beta$), we obtain $v\beta \in (X\beta)^* = X^*\beta$. Thus $v \in X^*$ and therefore $M = X^*$ as desired.
\qed

Now we can use Theorem \ref{fg} to derive a simplified description of $\rat(F)$:

\bt
\label{drs}
Let $L \subseteq F$. The following conditions are equivalent:
\bi
\item[(i)] $L \in {\rm Rat}(F)$;
\item[(ii)] $L$ is a finite union of subsets of the form
$$u_0X_1^*u_1X_2^*u_2\ldots X_n^*u_n$$
for $n \in \N$, $u_0,\ldots,u_n \in F$ and finite $X_1,\ldots,X_n \subset F$.
\ei
\et

\proof
Let $\cal{X}$ consist of all the subsets of $F$ satisfying condition (ii). It is immediate that ${\cal{X}} \subseteq \rat(F)$ and contains all the finite subsets of $F$. Thus it suffices to show that $\cal{X}$ is closed under (finite) union, product and the star operator.

Let $K,K' \in {\cal{X}}$. It is obvious that $K \cup K', KK' \in {\cal{X}}$. On the other hand, $K^* = X^*$ for some finite $X \subset F$ by Theorem \ref{fg}.
\qed

This is equivalent to say that $L \in \rat(F)$ if and only if $L = K\theta$ for some $K \in \rat(A^*)$ of {\em star height} 1 (see
\cite[Section I.6]{Sak}).

\section{Open problems}

The natural open problems in this context involve possible generalizations from the monogenic case to the finitely generated case, namely:

\bi
\item
are inclusion and equality decidable for rational subsets of $FIM_X$?
\item
is it decidable whether a rational subset of $FIM_X$ is recognizable?
\ei

With respect to Theorem \ref{fg}, it is easy to produce a counterexample:

\be
Let $X = \{ a,b\}$ and $L = (ab^*)^* \in {\rm Rat}(FIM_X)$. Then $L$ is not finitely generated.
\ee

Indeed, suppose that $L$ is finitely generated. Then $L \subseteq \{ a,ab,\ldots, ab^m\}^*$ for some $m \geq 1$. However, using Munn trees it is easy to see that $ab^{m+1} \notin \{ a,ab,\ldots, ab^m\}^*$, thus $L$ cannot be finitely generated.

\section*{Acknowledgement}

The author thanks the anonymous referee for his/her comments. The author was partially supported by CMUP, which is financed by national funds through FCT -- Funda\c c\~ao para a Ci\^encia e a Tecnologia, I.P., under the project with reference UIDB/00144/2020.

\end{document}